\documentclass[11pt,reqno]{amsart}
\usepackage{graphicx}
\usepackage{amsfonts}
\usepackage{mathrsfs}
 \usepackage{times}
 \usepackage{mathptmx}
\usepackage{amsmath}
\usepackage{txfonts}
\usepackage{amssymb}
\usepackage{amsmath,amsthm}
\usepackage{CJK,amsfonts}
\usepackage[body={7.5in,8.75in}, top=1.2in, right=0.9in,left=0.9in]{geometry}
\usepackage[colorlinks=true, linkcolor=red, citecolor=blue]{hyperref}
\allowdisplaybreaks[4] \numberwithin{equation}{section}

\newcommand{\N}{\mathbb N}
\newcommand{\R}{\mathbb R}

\newcommand{\abs}[1]{\left\lvert#1\right\rvert}
\newtheorem{thm}{Theorem}

\newtheorem{cor}[thm]{Corollary}
\newtheorem{lem}[thm]{Lemma}
\newtheorem{rem}[thm]{Remark}
\newtheorem{pro}[thm]{Proposition}

\DeclareMathOperator{\td}{d\mspace{-1mu}}
\def\squarebox#1{\hbox to #1{\hfill\vbox to #1{\vfill}}}

\begin{document}
\begin{CJK*}{GBK}{song}
\title[A note on fractional Askey--Wilson integrals]{A note on fractional Askey--Wilson integrals}
\author{Jian Cao$^{1}$ and Sama Arjika${}^{2}$ }
\dedicatory{\textsc}
\thanks{${}^1$Department
of Mathematics, Hangzhou Normal University, Hangzhou City, Zhejiang Province, 311121, China.   ${}^2$Department of Mathematics and Informatics, University of Agadez, Niger.}
\thanks{Email: 21caojian@hznu.edu.cn, rjksama2008@gmail.com.}



\begin{abstract}
In this paper, we generalize fractional $q$-integrals by the method of $q$-difference equation. In addition, we deduce fractional Askey--Wilson integral, reversal type fractional Askey--Wilson integral and Ramanujan type fractional Askey--Wilson integral.
\end{abstract}

\maketitle

\section{Introduction}

The fractional $q$-calculus theories have been applied successfully in many fields, which is a very suitable tool in describing and solving a lot of problems in numerous sciences \cite{podlubny,K-S-T}, such as high-energy physics, system control, biomedical engineering and economics. Its theoretical and applied research has become a hot spot in the world. The treatment from the point view of the $q$-calculus can open new perspectives as it did, for example, in optimal control problems \cite{Baleanu-Fernandez,bangerezako,Shahed-G-Al,Sayevand-T-B}. For more information, see details in \cite{cao-niu-jmaa17,Cao-HM-Liu,salam1966PEMS,eulerian,liu-16,wang23,wang08,wang09,wang22,mingjin,predrag-sladjana-miomir2007FCAA,predrag-sladjana-miomir2007AADM,predrag-sladjana-miomir2010AMS}.\par

In this paper, we follow the notations and terminology in
\cite{G-R} and suppose that $0<q<1$.
We first show a list of various definitions and notations in
$q$-calculus which are useful to understand the subject
of this paper. The basic hypergeometric series ${}_r\phi_s$ \cite{G-R}
\begin{equation}
{}_r\phi_s\biggl[\begin{matrix}
\begin{array}{ccc}
a_1,a_2,\ldots,a_r\\
b_1,b_2,\ldots,b_s
\end{array}
\end{matrix};q,z\biggr]=\sum_{n=0}^\infty\frac{\bigl(a_1,a_2,\ldots,a_r;q\bigr)_n}{\bigl(q,b_1,b_2,\ldots,b_s;q\bigr)_n}
\Bigl[(-1)^nq^{n\choose2}\Bigr]^{1+s-r}z^n
\end{equation}
converges absolutely for all $z$ if $r\le s$ and for $\abs{z}<1$ if $r=s+1$  and  for terminating. The compact factorials of ${}_r\phi_s$ are defined respectively by
\begin{equation}
(a;q)_0=1,\quad[a]_q:=\frac{1-q^a}{1-q},\quad (a;q)_n=\prod_{k=0}^{n-1}(1-aq^k),\quad
(a;q)_\infty=\prod_{k=0}^\infty(1-aq^k)
\end{equation}
and $(a_1,a_2,\dotsc,a_m;q)_n=(a_1;q)_n(a_2;q)_n\dotsm(a_m;q)_n$,
where $m\in\N:=\{1,2,3,\cdots\}\quad \text{and}\quad n\in\N_0:=\N\cup\{0\}$.
The $q$-gamma function is defined by \cite{G-R}
\begin{equation}
\Gamma_q(x)=\frac{(q;q)_\infty}{(q^x;q)_\infty}(1-q)^{1-x},\quad
x\in\R\backslash\{0,-1,-2,\ldots\}.
\end{equation}
The Thomae--Jackson $q$-integral is defined by \cite{G-R,jackson,thomae}\par
\begin{equation}\label{q-int-def}
\int_a^bf(x)\td_qx=(1-q)\sum_{n=0}^\infty\Bigl[bf(bq^n)-af(aq^n)\Bigr]q^n.
\end{equation}
The Riemann--Liouville fractional $q$-integral operator is introduced in \cite{salam66PEMS}
\begin{equation}\label{fra-def}
\Bigl(I_q^\alpha f\Bigr)(x)=\frac{x^{\alpha-1}}{\Gamma_q(\alpha)}
\int_0^x \bigl(qt/x;q\bigr)_{\alpha-1}f(t)\td_qt.
\end{equation}
The generalized Riemann--Liouville fractional $q$-integral operator is given by \cite{predrag-sladjana-miomir2007AADM}
\begin{equation}\label{fra-q-int}
\Bigl(I_{q,a}^\alpha f\Bigr)(x)=\frac{x^{\alpha-1}}{\Gamma_q(\alpha)}
\int_a^x \bigl(qt/x;q\bigr)_{\alpha-1}f(t)\td_qt,\quad\alpha\in \R^+.
\end{equation}
In fact, we rewrite fractional $q$-integral \eqref{fra-q-int} equivalently as follows by \eqref{q-int-def}
\begin{align}\label{fc-eq}
\Bigl(I_{q,a}^\alpha f\Bigr)(x)&=
\frac{x^{\alpha-1}(1-q)}{\Gamma_q(\alpha)}
\sum_{n=0}^\infty\Bigl[x\bigl(q^{n+1};q\bigr)_{\alpha-1}f\bigl(xq^n\bigr)
-a\bigl(aq^{n+1}/x;q\bigr)_{\alpha-1}f\bigl(aq^n\bigr)\Bigr]q^n.\nonumber
\end{align}

In paper \cite{cao2019}, author built the relations between the following fractional $q$-integrals and certain generating functions for $q$-polynomials.
\begin{pro}[{\cite[Theorem 3]{cao2019}}]\label{cao2019}
For $\alpha\in \R^+$ and $0<a<x<1$, if $\max\{\abs{at},\abs{az}\}<1$, we have
\begin{equation}\label{cao1}
\begin{aligned}
I_{q,a}^\alpha\biggl\{\frac{(bxz,xt;q)_\infty}{(xs,xz;q)_\infty}\biggr\}=\frac{(1-q)^\alpha(abz,at;q)_\infty}
{(as,az;q)_\infty}\sum_{k=0}^\infty
\frac{x^{\alpha+k}\bigl(a/x;q\bigr)_{\alpha+k}}{a^k(q;q)_{\alpha+k}}
{}_3\phi_2\left[\begin{matrix}
\begin{array}{cc}
q^{-k},as,az \\
at,abz\end{array}
\end{matrix};q,q\right].
\end{aligned}
\end{equation}
\end{pro}
In this paper, we generalize fractional $q$-integrals and give applications of  fractional Askey--Wilson integrals as follows.

\begin{thm}\label{cai1}
For $\alpha\in \R^+$ and $0<a<x<1$, if $\max\{\abs{at},\abs{az},\abs{aru}\}<1$, we have
\begin{equation}\label{cai1-1}
\begin{aligned}
I_{q,a}^\alpha\biggl\{\frac{(bxz,xt,xru;q)_\infty}{(xs,xz,xu;q)_\infty}\biggr\} =\frac{(1-q)^\alpha(abz,at,aru;q)_\infty}
{(as,az,au;q)_\infty}\sum_{k=0}^\infty
\frac{x^{\alpha+k}\bigl(a/x;q\bigr)_{\alpha+k}}{a^k(q;q)_{\alpha+k}}
{}_4\phi_3\left[\begin{matrix}
\begin{array}{cc}
q^{-k},as,az,au \\
abz,at,aru\end{array}
\end{matrix};q,q\right].
\end{aligned}
\end{equation}
\end{thm}
\begin{rem}
For $u=0$ in Theorem \ref{cai1} , equation \eqref{cai1-1} reduces to \eqref{cao1}.
\end{rem}\par
Before the proof of Theorem \ref{cai1}, the following lemmas are necessary.

\begin{lem}[{\cite[Proposition 1.2]{lu-09}}]
Let $f(a,b,c)$ be a three-variable analytic function in a neighborhood of $(a,b,c)=(0,0,0)\in C^3$. If $f(a,b,c)$ satisfies the difference equation
\begin{equation}\label{cai1-0}
(c-b)f(a,b,c)=abf(a,bq,cq)-bf(a,b,cq)+(c-ab)f(a,bq,c),
\end{equation}\par
then we have
\begin{equation}
f(a,b,c)=T(a,bD_c)\{f(a,0,c)\},
\end{equation}
where
\begin{equation}
T(a,bD_c)=\sum_{n=0}^\infty\frac{(a;q)_n}{(q;q)_n}(bD_c)^n,\quad
D_c\{f(c)\}=\frac{f(c)-f(cq)}{c}.
\end{equation}
\end{lem}

\begin{lem}\label{cai1-2}
For $\max\{\abs{as},\abs{az},\abs{au}\}<1$, we have
\begin{multline}\label{cai1-3}
(s-u)\frac{(abz,at,aru;q)_\infty}{(as,az,au)_\infty} 
=ur\frac{(abz,at,aruq;q)_\infty}{(asq,az,auq)_\infty}
-u\frac{(abz,at,aru;q)_\infty}{(asq,az,au)_\infty}
+(s-ur)\frac{(abz,at,aruq;q)_\infty}{(as,az,auq)_\infty}.
\end{multline}
\end{lem}

\textbf{Proof of Lemma \ref{cai1-2}.}
The right-hand side (RHS) of equation \eqref{cai1-3} equals
\begin{multline}
ur\frac{(1-xs)(1-xu)}{1-xur}\frac{(abz,at,aru;q)_\infty}{(as,az,au)_\infty}
-u(1-xs)\frac{(abz,at,aru;q)_\infty}{(as,az,au)_\infty}\\+(s-ur)\frac{1-xu}{1-xur}\frac{(abz,at,aru;q)_\infty}{(as,az,au)_\infty}
=(s-u)\frac{(abz,at,aru;q)_\infty}{(as,az,au)_\infty},
\end{multline}
which is equals to the left-hand side (LHS) of the equation \eqref{cai1-3}. The proof is complete.

\begin{lem}[{\cite[Eq.~(2.3)]{chen-gu}}]\label{T01}
For $\max\{\abs{bt},\abs{ct}\}<1$, we have
\begin{equation}
T(a,bD_c)\left\{\frac{1}{(ct;q)_\infty}\right\}=\frac{(abt;q)_\infty}{(bt,ct;q)_\infty}.
\end{equation}
\end{lem}
\begin{proof}[Proof of Theorem \ref{cai1}]
Denoting the RHS of the equation \eqref{cai1-1} by $f(r,u,s)$, and rewriting $f(r,u,s)$ equivalently by
\begin{align*}
f(r,u,s)=\sum_{k=0}^\infty
\frac{x^{\alpha+k}\bigl(a/x;q\bigr)_{\alpha+k}}{a^k(q;q)_{\alpha+k}}
\sum_{j=0}^k\frac{(q^{-k};q)_jq^j}{(q;q)_j}\cdot\frac{(1-q)^\alpha(abzq^j,atq^j,aruq^j;q)_\infty}
{(asq^j,azq^j,auq^j;q)_\infty},
\end{align*}
we check that $f(r,u,s)$ satisfies the equation \eqref{cai1-0} by Lemma \ref{cai1-2}, then we have
\begin{align*}
f(r,u,s)&=T(r,uD_s)\left\{f(r,0,s)\right\}\\
&=T(r,uD_s)\left\{\sum_{k=0}^\infty
\frac{x^{\alpha+k}\bigl(a/x;q\bigr)_{\alpha+k}}{a^k(q;q)_{\alpha+k}}
\sum_{j=0}^k\frac{(q^{-k};q)_jq^j}{(q;q)_j}\cdot\frac{(1-q)^\alpha(abzq^j,atq^j;q)_\infty}
{(asq^j,azq^j;q)_\infty}\right\}\\
&=T(r,uD_s)\left\{I_{q,a}^\alpha\Biggl\{\frac{(xbz,xt;q)_\infty}
{(xs,xz;q)_\infty}\Biggr\}\right\}\\
&=I_{q,a}^\alpha \left\{T(r,uD_s)\Biggl\{\frac{(xbz,xt;q)_\infty}
{(xs,xz;q)_\infty}\Biggr\}\right\}\\
&=I_{q,a}^\alpha\left\{\frac{(xbz,xt;q)_\infty}
{(xz;q)_\infty}\cdot T(r,uD_s)\Biggl\{\frac{1}{(xs;q)_\infty}\Biggr\}\right\},
\end{align*}
which equals the LHS of the equation \eqref{cai1-1} by Lemma \ref{T01}. The proof is complete.
\end{proof}
The rest of the paper is organized as follows.
In Section \ref{section2}, we give the fractional Askey--Wilson integral. In Section \ref{section3}, we obtain the reversal type fractional Askey--Wilson integral. In Section \ref{section4}, we get the Ramanujan type fractional Askey--Wilson integral.

\section{A generalization of Askey--Wilson integrals}\label{section2}
In 1985, Askey and Wilson gave the famous Askey--Wilson integral, which greatly promoted the research and development of orthogonal polynomials. Chen and Gu \cite{chen-gu}, Liu \cite{liu-09} and Cao \cite{cao-jmaa14} et al have promoted Askey--Wilson integral by different methods. For more information, see details in \cite{A-W-85,cao-jmaa14,chen-gu,liu-09}. \par
In this section, we use the fractional $q$-integrals to generalize Askey--Wilson integrals.
\begin{pro}[{\cite[Theorem 2.1]{A-W-85}}]
If max$\{\abs{a},\abs{b},\abs{c},\abs{d}\}<1$, we have
\begin{equation}\label{A-W-I}
\int_0^\pi\frac{h(\cos2\theta;1)}{h(\cos\theta;a,b,c,d)}d\theta
=\frac{2\pi(abcd;q)_\infty}{(q,ab,ac,ad,bc,bd,cd;q)_\infty},
\end{equation}
where
\begin{align*}
h(\cos\theta;a)&=(ae^{i\theta},ae^{-i\theta};q)_\infty,\\
h(\cos\theta;a_1,a_2, \cdots ,a_m)&=h(\cos\theta;a_1)h(\cos\theta;a_2)\cdots h(\cos\theta;a_m).
\end{align*}
\end{pro}

\begin{thm}\label{cai2-0}For $\alpha\in \R^+$, if $\max\{\abs{a},\abs{b},\abs{c},\abs{d}\}<1$, we have
\begin{multline}\label{cai2-1}
\int_0^\pi\frac{h(\cos2\theta;1)}{h(\cos\theta;a,b,c,d)}
\sum_{k=0}^\infty
\frac{x^{\alpha+k}\bigl(a/x;q\bigr)_{\alpha+k}}{a^k(q;q)_{\alpha+k}}
{}_4\phi_3\left[\begin{matrix}
\begin{array}{cc}
q^{-k},abcd,ae^{i\theta},ae^{-i\theta} \\
ab,ac,ad\end{array}
\end{matrix};q,q\right]d\theta\\
=\frac{2\pi(abcd;q)_\infty}{(q,ab,ac,ad,bc,bd,cd;q)_\infty}
\frac{x^\alpha(a/x;q)_\alpha}{(q;q)_\alpha}.
\end{multline}
\end{thm}

\begin{cor}For $\alpha\in \R^+$, if $max\{\abs{a},\abs{b},\abs{c}\}<1$, we have
\begin{multline}\label{cai2-2}
\int_0^\pi\frac{h(\cos2\theta;1)}{h(\cos\theta;a,b,c)}
\sum_{k=0}^\infty
\frac{x^{\alpha+k}\bigl(a/x;q\bigr)_{\alpha+k}}{a^k(q;q)_{\alpha+k}}
{}_3\phi_2\left[\begin{matrix}
\begin{array}{cc}
q^{-k},ae^{i\theta},ae^{-i\theta} \\
ab,ac\end{array}
\end{matrix};q,q\right]d\theta
=\frac{2\pi}{(q,ab,ac,bc;q)_\infty}
\frac{x^\alpha(a/x;q)_\alpha}{(q;q)_\alpha}.
\end{multline}
\end{cor}
\begin{rem}
For $d=0$ in Theorem \ref{cai2-0}, equation \eqref{cai2-1} reduces to \eqref{cai2-2}.
\end{rem}
\begin{proof}[Proof of Theorem \ref{cai2-0}]
The equation \eqref{cai2-1} can be rewrite equivalently by
\begin{align}\label{cai2-3}
&\int_0^\pi\frac{h(\cos2\theta;1)}{h(\cos\theta;b,c,d)}
\frac{(ab,ac,ad;q)_\infty}{(ae^{i\theta},ae^{-i\theta},abcd;q)_\infty}
\sum_{k=0}^\infty
\frac{x^{\alpha+k}\bigl(a/x;q\bigr)_{\alpha+k}}{a^k(q;q)_{\alpha+k}}\nonumber\\
&\quad\quad\quad\quad\times{}_4\phi_3\left[\begin{matrix}
\begin{array}{cc}
q^{-k},abcd,ae^{i\theta},ae^{-i\theta} \\
ab,ac,ad\end{array}
\end{matrix};q,q\right]d\theta
=\frac{2\pi}{(q,bc,bd,cd;q)_\infty}
\frac{x^\alpha(a/x;q)_\alpha}{(q;q)_\alpha}.
\end{align}
By using Theorem \ref{cai1}, the LHS of the equation \eqref{cai2-3} is equivalent to
\begin{align*}
\int_0^\pi\frac{h(\cos2\theta;1)}{h(\cos\theta;b,c,d)}
I_{q,a}^\alpha\biggl\{\frac{(xb,xc,xd;q)_\infty}
{(xe^{i\theta},xe^{-i\theta},xbcd;q)_\infty}\biggr\}d\theta& =I_{q,a}^\alpha\biggl\{\int_0^\pi\frac{h(\cos2\theta;1)(xb,xc,xd;q)_\infty}
{h(\cos\theta;x,b,c,d)(xbcd;q)_\infty}d\theta \biggr\}\\
&\quad\quad =I_{q,a}^\alpha\biggl\{\frac{2\pi}{(q,bc,bd,cd;q)_\infty}\biggr\}\\
&\quad\quad =\frac{2\pi}{(q,bc,bd,cd;q)_\infty}
\frac{x^\alpha(a/x;q)_\alpha}{(q;q)_\alpha},
\end{align*}
which is the RHS of the equation \eqref{cai2-3}. The proof is complete.

\end{proof}
\section{A generalization of reversal type Askey--Wilson integrals}\label{section3}
In this section, we use the fractional $q$-integrals to expand reversal Askey--Wilson integrals.
\begin{pro}[Reversal Askey--Wilson integral {\cite{R-A-W-87}}]
For $\abs{qabcd}<1$, there holds
\begin{align}\label{R-A-W}
\int_{-\infty}^\infty\frac{h(i \sinh x;qa,qb,qc,qd)}{h(\cosh2x;-q)}dx
=\frac{(q,qab,qac,qad,qbc,qbd,qcd;q)_\infty}{(qabcd;q)_\infty}\log(q^{-1}),
\end{align}
where
\begin{equation}
h(i \sinh \alpha x;t)=\prod_{k=0}^\infty(1-2iq^kt\sinh \alpha x+q^{2k}t^2)=(ite^{\alpha x},-ite^{-\alpha x};q)_\infty.\nonumber
\end{equation}
\end{pro}

\begin{thm}\label{cai3-0}
For $\alpha\in \R^+$ and $ \abs{qabcd}<1$, we have
\begin{align}\label{cai3-1}
&\int_{-\infty}^\infty\frac{h(i \sinh t;qa,qb,qc,qd)}{h(\cosh2t;-q)}
\sum_{k=0}^\infty
\frac{x^{\alpha+k}\bigl(a/x;q\bigr)_{\alpha+k}}{a^k(q;q)_{\alpha+k}}
{}_4\phi_3\left[\begin{matrix}
\begin{array}{cc}
q^{-k},qab,qac,qad \\
iaqe^t,-iaqe^{-t},qabcd\end{array}
\end{matrix};q,q\right] dt\nonumber\\
&\quad=\frac{(q,qab,qac,qad,qbc,qbd,qcd;q)_\infty}{(qabcd;q)_\infty}
\frac{x^\alpha(a/x;q)_\alpha}{(q;q)_\alpha}\log(q^{-1}).
\end{align}
\end{thm}

\begin{cor}
For $\alpha\in \R^+$, we have
\begin{align}\label{cai3-2}
\int_{-\infty}^\infty\frac{h(i \sinh t;qa,qb,qc)}{h(\cosh2t;-q)}
\sum_{k=0}^\infty
\frac{x^{\alpha+k}\bigl(a/x;q\bigr)_{\alpha+k}}{a^k(q;q)_{\alpha+k}}
{}_3\phi_2\left[\begin{matrix}
\begin{array}{cc}
q^{-k},qab,qac \\
iaqe^t,-iaqe^{-t}\end{array}
\end{matrix};q,q\right] dt\nonumber\\
=\frac{(q,qab,qac,qbc;q)_\infty x^\alpha(a/x;q)_\alpha}{(q;q)_\alpha}\log(q^{-1}).
\end{align}
\end{cor}

\begin{rem}
For $d=0$ in Theorem \ref{cai3-0}, equation \eqref{cai3-1} reduces to \eqref{R-A-W}.
\end{rem}
\begin{proof}[Proof of Theorem \ref{cai3-0}]
The equation \eqref{cai3-1} can be rewrite
\begin{multline}\label{cai3-2}
\int_{-\infty}^\infty\frac{h(i \sinh t;qb,qc,qd)}{h(\cosh2t;-q)}\frac{(iaqe^t,-iaqe^{-t},qabcd;q)_\infty}{(qab,qac,qad;q)_\infty}
\sum_{k=0}^\infty
\frac{x^{\alpha+k}\bigl(a/x;q\bigr)_{\alpha+k}}{a^k(q;q)_{\alpha+k}} {}_4\phi_3\left[\begin{matrix}
\begin{array}{cc}
q^{-k},qab,qac,qad \\
iaqe^t,-iaqe^{-t},qabcd\end{array}
\end{matrix};q,q\right] dt\\
=(q,qbc,qbd,qcd;q)_\infty
\frac{x^\alpha(a/x;q)_\alpha}{(q;q)_\alpha}\log(q^{-1}),
\end{multline}
By using Theorem \ref{cai1}, the LHS of the equation \eqref{cai3-2} equals
\begin{align*}
\int_{-\infty}^\infty\frac{h(i \sinh t;qb,qc,qd)}{h(\cosh2t;-q)}
I_{q,a}^\alpha\biggl\{\frac{(ixqe^t,-ixqe^{-t},qxbcd;q)_\infty}
{(qxb,qxc,qxd;q)_\infty}\biggr\}dt&=I_{q,a}^\alpha\biggl\{\int_{-\infty}^\infty\frac{h(i \sinh t;qx,qb,qc,qd)(qxbcd;q)_\infty}{h(\cosh2t;-q)(qxb,qxc,qxd;q)_\infty}dt\biggr\}\\
&\quad\quad\quad\quad=I_{q,a}^\alpha\biggl\{(q,qbc,qbd,qcd;q)_\infty\log(q^{-1})\biggr\}\\
&\quad\quad\quad\quad=(q,qbc,qbd,qcd;q)_\infty
\frac{x^\alpha(a/x;q)_\alpha}{(q;q)_\alpha}\log(q^{-1}).
\end{align*}
which is the RHS of the equation \eqref{cai3-2}. The proof is complete.

\end{proof}

\section{A generalization of Ramanujan type Askey--Wilson integrals }\label{section4}
In 1994, Atakishiyev discovered the Atakishiyev integral by Ramanujan's method. After that, Wang \cite{mingjin} and Liu \cite{liu-15} et al conducted different generalization and proving lines for this integral.\par
In this part, we get the generating functions of Atakishiyev integral by fractional $q$-integral.
\begin{pro}[Atakishiyev integral {\cite{A-I-94}}]
If $\alpha$ is a real number and $q=exp(-2\alpha^2)$, then we have
\begin{align}\label{R-A-I}
\int_{-\infty}^\infty h(i \sinh \alpha t;a,b,c,d)e^{-t^2}\cosh \alpha t dt
=\sqrt\pi q^{-\frac{1}{8}}\frac{(ab/q,ac/q,ad/q,bc/q,bd/q,cd/q;q)_\infty}{(abcd/q^3;q)_\infty}.
\end{align}
\end{pro}

\begin{thm}\label{cai4-0}
For $\alpha\in \R^+$ and $\abs{abcd/q^3}<1$, if $\alpha$ is a real number and $q=exp(-2\alpha^2)$, then we have
\begin{align}\label{cai4-1}
&\int_{-\infty}^\infty h(i \sinh \alpha t;a,b,c,d)e^{-x^2}\cosh \alpha t
\sum_{k=0}^\infty
\frac{x^{\alpha+k}\bigl(a/x;q\bigr)_{\alpha+k}}{a^k(q;q)_{\alpha+k}} {}_4\phi_3\left[\begin{matrix}
\begin{array}{cc}
q^{-k},ab/q,ac/q,ad/q \\
iae^{\alpha t},-iae^{-\alpha t},abcd/q^3\end{array}
\end{matrix};q,q\right]dt \nonumber \\
&\qquad=\sqrt\pi q^{-\frac{1}{8}}\frac{(ab/q,ac/q,ad/q,bc/q,bd/q,cd/q;q)_\infty}{(abcd/q^3;q)_\infty}
\frac{x^\alpha(a/x;q)_\alpha}{(q;q)_\alpha}.
\end{align}
\end{thm}

\begin{cor}
For $\alpha\in \R^+$, if $\alpha$ is a real number and $q=exp(-2\alpha^2)$, then we have
\begin{align}\label{cai4-6}
\int_{-\infty}^\infty h(i \sinh \alpha t;a,b,c)&e^{-x^2}\cosh \alpha t
\sum_{k=0}^\infty
\frac{x^{\alpha+k}\bigl(a/x;q\bigr)_{\alpha+k}}{a^k(q;q)_{\alpha+k}}
{}_3\phi_2\left[\begin{matrix}
\begin{array}{cc}
q^{-k},ab/q,ac/q \\
iae^{\alpha t},-iae^{-\alpha t}\end{array}
\end{matrix};q,q\right]dt \nonumber \\
&=\sqrt\pi q^{-\frac{1}{8}}(ab/q,ac/q,bc/q;q)_\infty
\frac{x^\alpha(a/x;q)_\alpha}{(q;q)_\alpha}.
\end{align}
\end{cor}

\begin{rem}
For $d=0$ in Theorem \ref{cai4-0}, equation \eqref{cai4-1} reduces to \eqref{cai4-6}.
\end{rem}

\begin{proof}[Proof of Theorem \ref{cai4-0}]
The equation \eqref{cai4-1} can be rewrited
\begin{multline}\label{cai4-2}
\int_{-\infty}^\infty h(i \sinh \alpha t;b,c,d)
e^{-x^2}\cosh \alpha t\frac{(iae^{\alpha t},-iae^{-\alpha t},abcd/q^3;q)_\infty}{(ab/q,ac/q,ad/q;q)_\infty}
\sum_{k=0}^\infty
\frac{x^{\alpha+k}\bigl(a/x;q\bigr)_{\alpha+k}}{a^k(q;q)_{\alpha+k}}\\
\times{}_4\phi_3\left[\begin{matrix}
\begin{array}{cc}
q^{-k},ab/q,ac/q,ad/q \\
iae^{\alpha t},-iae^{-\alpha t},abcd/q^3\end{array}
\end{matrix};q,q\right] dt\\
=\sqrt\pi q^{-\frac{1}{8}}(bc/q,bd/q,cd/q;q)_\infty
\frac{x^\alpha(a/x;q)_\alpha}{(q;q)_\alpha}.
\end{multline}
By using Theorem \ref{cai1}, the LHS of the equation \eqref{cai4-2} is equal to
\begin{align*}
&\int_{-\infty}^\infty h(i \sinh \alpha t;b,c,d)
e^{-x^2}\cosh \alpha t\cdot
I_{q,a}^\alpha\biggl\{\frac{(ixe^{\alpha t},-ixe^{-\alpha t},xbcd/q^3;q)_\infty}{(xb/q,xc/q,xd/q;q)_\infty}\biggr\}dt\\
&\quad\quad\quad\quad=I_{q,a}^\alpha\biggl\{\int_{-\infty}^\infty h(i \sinh \alpha t;x,b,c,d)
e^{-x^2}\cosh \alpha t\cdot\frac{(xbcd/q^3;q)_\infty}{(xb/q,xc/q,xd/q;q)_\infty}dt\biggr\}\\
&\quad\quad\quad\quad=I_{q,a}^\alpha\biggl\{\sqrt\pi q^{-\frac{1}{8}}(bc/q,bd/q,cd/q;q)_\infty\biggr\}\\
&\quad\quad\quad\quad=\sqrt\pi q^{-\frac{1}{8}}(bc/q,bd/q,cd/q;q)_\infty
\frac{x^\alpha(a/x;q)_\alpha}{(q;q)_\alpha},
\end{align*}
which is the RHS of the equation \eqref{cai4-2}. The proof is complete.
\end{proof}
\section*{Acknowledgments}
The author would like to thank the referees and editors for their many valuable comments and suggestions. This work was supported by the Zhejiang Provincial Natural Science Foundation of China (No.~LY21A010019).

\end{CJK*}
\end{document}